\documentclass[11pt,a4paper]{article}

\usepackage{amssymb}

\usepackage{epsfig}
\usepackage{url}

\newenvironment{keywords}{ \noindent {\small\bf Key Words}:}{ }

\def\be{\begin{equation}}
\def\ee{\end{equation}}

\newtheorem{thm}{Theorem}

\newtheorem{dfn}{Definition}

\newtheorem{pf}{Proof}

\begin{document}

\title{Lipschitz gradients for global optimization\\ in a
one-point-based partitioning scheme\footnote{This work was supported
by the grants 1960.2012.9 and MK-3473.2010.1 awarded by the
President of the Russian Federation for supporting the leading
research groups and young researchers, respectively, as well as by
the grant 11-01-00682-a awarded by the Russian Foundation for
Fundamental Research. The authors thank also the Italian
Inter-University Consortium for the Application of Super-Computing
for Universities and Research (CASPUR), project ``Advanced
techniques for global optimization: Numerical methods and
applications'' in the framework of the ``HPC Grant 2011 on GPU
cluster''.}}

\newcommand{\nms}{\normalsize}
\author{Dmitri E. Kvasov\footnote{kvadim@si.deis.unical.it} and  Yaroslav D. Sergeyev\footnote{Corresponding author,
yaro@si.deis.unical.it  }
\\ \\ [-2pt]
      \nms DEIS -- University of Calabria,\\[-4pt]
       \nms   Via P. Bucci 42C, 87036 -- Rende (CS), Italy,\\[-4pt]
       \nms       and  \\[-4pt]
       \nms   Software Department, N.I.\,Lobachevsky,\\[-4pt]
       \nms   State University, Nizhni Novgorod, Russia
}

\date{}

\maketitle

\begin{abstract}
A global optimization problem is studied where the objective
function~$f(x)$ is a multidimensional black-box function and its
gradient $f'(x)$ satisfies the Lipschitz condition over a
hyperinterval with an unknown Lipschitz constant~$K$. Different
methods for solving this problem by using an a priori given
estimate of $K$, its adaptive estimates, and adaptive estimates of
local Lipschitz constants are known in the literature. Recently,
the authors have proposed  a one-dimensional algorithm working
with multiple estimates of the Lipschitz constant for $f'(x)$ (the
existence of such an algorithm was a challenge for 15 years).  In
this paper, a new multidimensional geometric method evolving the
ideas of this one-dimensional scheme and using an efficient
one-point-based partitioning strategy is proposed. Numerical
experiments executed on  800 multidimensional test functions
demonstrate quite a promising performance in comparison with
popular DIRECT-based methods.
\end{abstract}


 \begin{flushleft}
\textbf{MSC classes:}  65K05, 90C26, 90C56.
\end{flushleft}

\section{Introduction}

Global optimization is an important part of numerical analysis
(see, e.g., \cite{Pardalos:et:al.(2000), Strongin&Sergeyev(2000),
Trigiante(2000), Zhigljavsky&Zilinskas(2008)}). It considers
problems derived from complex industrial applications where the
objective function $f(x)$ to be minimized is defined over a
hyperinterval $D \subset R^N$, $N \geq 1$, and can be black-box,
multiextremal, and requiring high computational resources for its
evaluation (see, e.g.,~\cite{Audet:et:al.(2005),
Dumas:et:al.(2009), Luo:et:al.(2011), Mockus(2000),
Pardalos:et:al.(2000), Pinter(1996), Sergeyev&Kvasov(2008),
Strongin&Sergeyev(2000), Zhigljavsky&Zilinskas(2008)}). Solving
efficiently this type of problems is a great challenge, since they
present a high number of local minimizers (only a few of which can
be global ones), often with extremely different values, and do not
present a simple mathematical description of the global optima.

One of the natural and powerful (from both the theoretical and the
applied points of view) assumptions on these problems is that the
objective function has bounded slopes, i.e.,
 \be \label{LGOP_L}
   | f(x') - f(x'') | \le L \| x'- x'' \|,
   \hspace*{3mm}   x', x''\in D, \hspace*{3mm} 0 <L < \infty,
 \ee
where $\|\cdot\|$ denotes, usually, the Euclidean norm (other
norms can be also used, see, e.g., \cite{Gaviano&Lera(2008),
Nesterov(2004)}) and $L$ is the (unknown) Lipschitz constant. In
this case,   Lipschitz global optimization  methods can be applied
(see, e.g., \cite{Evtushenko(1985), Horst&Pardalos(1995),
Horst&Tuy(1996), Pinter(1996), Sergeyev&Kvasov(2008),
Strongin&Sergeyev(2000), Zhigljavsky&Zilinskas(2008)} and the
references given therein). They can be distinguished, for example,
by the way in which information about the Lipschitz constant is
obtained and by the strategy of exploration of the admissible
region.

In the literature, there exist at least four approaches to specify
the Lipschitz constant $L$ from~(\ref{LGOP_L}): (i) it can be
given a priori (see, e.g.,~\cite{Evtushenko:et:al.(2009),
Evtushenko(1985), Horst&Pardalos(1995)}); (ii)~its global estimate
over the whole domain can be used (see,
e.g.,~\cite{Horst&Tuy(1996), Kvasov&Sergeyev(2003), Pinter(1996),
Strongin&Sergeyev(2000)}); (iii)~local Lipschitz constants can be
estimated (see, e.g.,~\cite{Kvasov:et:al.(2003),
Sergeyev&Kvasov(2008), Sergeyev(1995), Strongin&Sergeyev(2000)});
(iv)~several estimates of $L$ can be chosen from a set of possible
values (see, e.g.,~\cite{Finkel&Kelley(2006),
Gablonsky&Kelley(2001), He:et:al.(2002), Jones:et:al.(1993),
Liuzzi:et:al.(2010), Sergeyev&Kvasov(2006),
Sergeyev&Kvasov(2008)}). In their work global optimization methods
using multiple estimates of the Lipschitz constants have proved to
be particularly attractive for studying applied problems (see,
e.g.,~\cite{DiSerafino:et:al.(2011), Graf:et:al.(2007),
He:et:al.(2002), Moles:et:al.(2003), Panning:et:al.(2008)}, other
references can be found, e.g., in~\cite{Kvasov&Sergeyev(2009),
Sergeyev&Kvasov(2006)}).

In exploring the multidimensional search domain, various adaptive
partitioning strategies can be applied. For example,
one-point-based algorithms subsequently subdivide the search
region in smaller ones and evaluate the objective function at one
point within each subregion (see,
e.g.,~\cite{Evtushenko&Posypkin(2011), Evtushenko(1985),
Gablonsky&Kelley(2001), Jones:et:al.(1993), Sergeyev(2005)}).
Partitions of the search domain into hyperintervals, based on
evaluating the objective function at the two vertices
corresponding to the main diagonal of hyperintervals called
\textit{diagonal partitioning strategies}, can also be
successfully used  (see, e.g.,\cite{Gergel(1997),
Kvasov:et:al.(2003), Kvasov&Sergeyev(2003), Pinter(1996),
Sergeyev&Kvasov(2006), Sergeyev&Kvasov(2008)}). More complex
partitions, based on simplices, auxiliary functions of various
nature, and so on, have also been proposed (see,
e.g.,~\cite{Horst&Tuy(1996), Lera&Sergeyev(2010b),
Wu:et:al.(2005), Zhigljavsky&Zilinskas(2008)}; many other
references can be found in~\cite{Sergeyev&Kvasov(2011)}).

The choice of the regions to be partitioned is based on an
information about the objective function obtained during the
search. It can be either of the probabilistic type (e.g., Bayesian
approach applying the theory of random functions to a mathematical
representation of available (certain or uncertain) a priori
information on the objective function behavior, see,
e.g.,~\cite{Jones:et:al.(1998), Lera&Sergeyev(2010a),
Mockus(2000), Strongin&Sergeyev(2000),
Zhigljavsky&Zilinskas(2008)}), or of the deterministic one (e.g.,
geometric approach making a use of different auxiliary functions
to estimate the behavior of $f(x)$ over the search region, see,
e.g., \cite{Baritompa(1993), Horst&Pardalos(1995),
Horst&Tuy(1996), Pinter(1996), Sergeyev&Kvasov(2008),
Sergeyev(1995), Sergeyev(1998), Strongin&Sergeyev(2000),
Zhigljavsky&Zilinskas(2008)}).

In this paper, a particular class of the Lipschitz global
optimization problems is considered, namely, the class of problems
with differentiable objective functions having the Lipschitz
gradients $f'(x)$, i.e.,

 \be \label{LGOP_f}
   f ^{*}=f(x^{*})=\min_{x \in D} \; f(x),
 \ee
 \be \label{LGOP_K}
   \| f'(x') - f'(x'') \| \le K \| x'- x'' \|,
   \hspace*{3mm}   x', x''\in D, \hspace*{3mm} 0 <K < \infty,
 \ee
where
 \be \label{LGOP_D}
   D = [a,b] = \{ x \in R^N: a(j) \leq x(j) \leq b(j) \}.
 \ee

It is supposed in this formulation that the objective function
$f(x)$ can be black-box, multiextremal, its gradient
$f'(x)=\left(\frac{\partial f(x)}{\partial x(1)}, \frac{\partial
f(x)}{\partial x(2)}, \ldots, \frac{\partial f(x)}{\partial
x(N)}\right)^{T}$ (which could be itself a costly multiextremal
black-box vector-function) can be calculated during the search,
and $f'(x)$ is Lipschitz-continuous with some fixed, but unknown,
constant~$K$, $0<K<\infty$, over $D$. These problems are often
encountered in engineering applications (see, e.g.,
\cite{Pinter(1996), Sergeyev&Kvasov(2008),
Strongin&Sergeyev(2000)}), particularly, in electrical engineering
optimization problems (see, e.g., \cite{Sergeyev:et:al.(1999),
Sergeyev&Kvasov(2008), Strongin&Sergeyev(2000)}).

In the literature, several methods for solving this problem have
been proposed. They can be also distinguished, for instance, with
respect to the way the Lipschitz constant~$K$ is estimated in
their work. There exist algorithms using an a priori given
estimate of $K$ (see, e.g., \cite{Baritompa(1993),
Breiman&Cutler(1993), Sergeyev(1998)}), its adaptive estimates
(see, e.g., \cite{Gergel(1997), Sergeyev&Kvasov(2008),
Sergeyev(1998)}), and adaptive estimates of local Lipschitz
constants (see, e.g., \cite{Sergeyev&Kvasov(2008),
Sergeyev(1998)}). Algorithms working with a number of Lipschitz
constants for~$f'(x)$  chosen from a set of possible values
varying from zero to infinity were not known till 2009 when such
an algorithm for solving the one-dimensional problem
(\ref{LGOP_f})--(\ref{LGOP_D}) has been proposed
in~\cite{Kvasov&Sergeyev(2009)}. Its extension to the
multidimensional case is not a trivial task in contrast to the
DIRECT~method (see~\cite{Jones:et:al.(1993)}) proposed in 1993 for
solving problems with the Lipschitz objective function.

The present paper solves this more than 15-year open problem of
constructing multidimensional global optimization methods working
with multiple estimates of the Lipschitz constants for~$f'(x)$. A
new multidimensional geometric method for finding solutions to the
problem (\ref{LGOP_f})--(\ref{LGOP_D}) is introduced and studied
here. It uses a new one-point-based partitioning strategy
(see~\cite{Sergeyev&Kvasov(2008), Sergeyev(2005)}) and works with
a number of estimates of the Lipschitz constant $K$ for $f'(x)$.
Such multiple (from zero to infinity) estimates of $K$
from~(\ref{LGOP_K}) are used to calculate the lower bounds of the
objective function over the hyperintervals of a current partition
of the search domain and to produce new trial points (i.e., points
at which both the objective function~$f(x)$ and its
gradient~$f'(x)$ are evaluated). In the framework of geometric
algorithms, this kind of estimating the Lipschitz constant can be
interpreted as examination of all admissible minorant functions
during the current iteration of the algorithm without constructing
a specific one. A particular attention in the new algorithm is
given to the improvement of the current minimal function value
(the so-called \textit{record value}) in order to provide a faster
convergence to a global minimizer. As demonstrated by extensive
numerical experiments executed on 800 test functions from the
differentiable GKLS test classes
(see~\cite{Gaviano:et:al.(2003)}), the usage of gradients allows
one to obtain, as expected, an acceleration in comparison with the
DIRECT-based methods.

The paper is organized as follows. In Section~2, a theoretical
background of the new algorithm is presented. Section~3 is
dedicated to the description of the algorithm and to its
convergence analysis. Finally, Section~4 contains results of
numerical experiments executed on 800 test functions.

\section{Theoretical background}

In this section, the main theoretical results, necessary for
introducing the new algorithm, are obtained. First, a new
partitioning strategy developed in the framework of the
one-point-based partition approach is described. The second part
presents a technique for estimating the lower bounds of the
objective function over hyperintervals. The third part is
dedicated to the introduction of a procedure for determining
nondominated hyperintervals, i.e., hyperintervals having the
smallest lower bound for some particular estimate of the Lipschitz
constant for $f'(x)$. They are candidates for partitioning at each
iteration of the new method.

\subsection{One-point-based partitioning strategy} \label{sectionNewStrategy}

In this section, a new efficient one-point-based partitioning
scheme proposed in~\cite{Sergeyev(2005)} (see
also~\cite{Sergeyev&Kvasov(2008)}) is considered which is based on
a diagonal partitioning strategy from~\cite{Sergeyev&Kvasov(2008),
Sergeyev(2000)}. In this scheme, the function $f(x)$ and its
gradient $f'(x)$ are evaluated only at one vertex (either $a_i$ or
$b_i$) of the main diagonal of each hyperinterval~$D_i=[a_i, b_i]$
of the current partition independently of the problem dimension
(recall that performing each trial is a time-consuming operation).

\begin{figure}[t]
\centerline{\epsfig{file=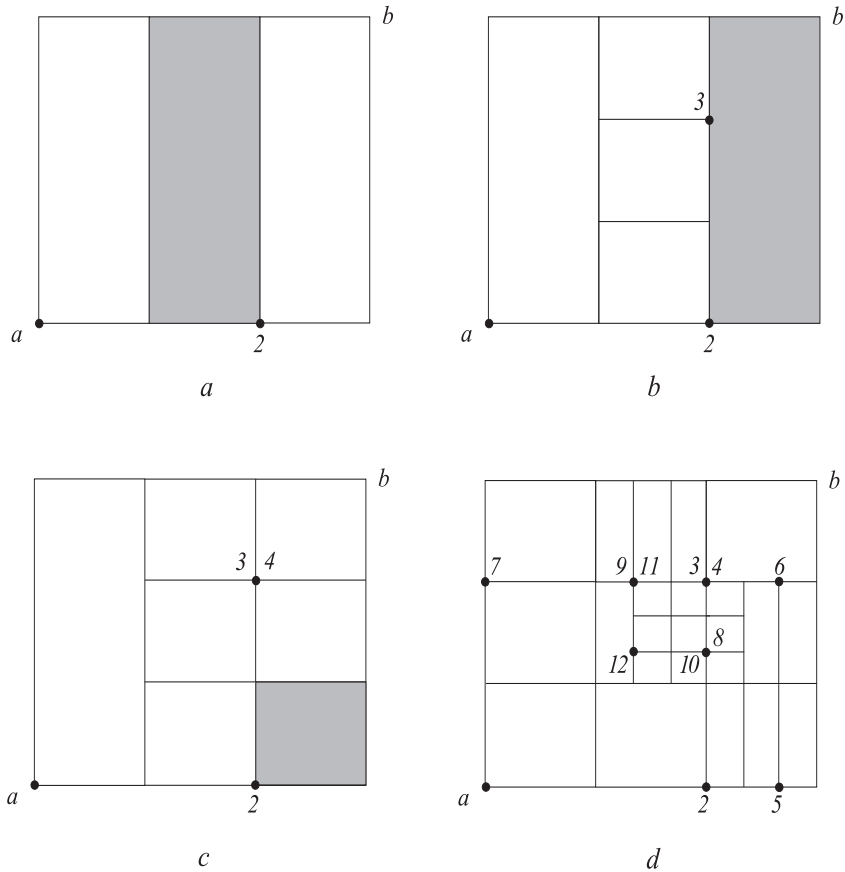,width=110mm,height=110mm,angle=0,silent=}}
\caption{An example of subdivisions by a new one-point-based
partitioning strategy} \label{fig_OnePoint}
\end{figure}

Let us start the description of this scheme with a two-dimensional
example shown in Fig.~\ref{fig_OnePoint}. In this Figure,
partitions of the admissible region~$D$ produced by the algorithm
at several initial iterations are presented starting from the
first trial at the point $a$ (it is supposed here that a single
iteration consists of the subdivision of only one hyperinterval).
Black dots represent the trial points and the numbers around these
dots indicate iterations at which these trial points have been
generated. The terms `interval' and `subinterval' will be used to
denote two-dimensional rectangular domains.

In Fig.~\ref{fig_OnePoint}a, the situation after the first two
iterations is presented. Particularly, at the second iteration,
the interval $D$ is partitioned into three subintervals of equal
area (equal volume in a general case). This subdivision is
performed by two lines (hyperplanes) orthogonal to the longest
edge of~$D$ (see Fig.~\ref{fig_OnePoint}a). The trial (evaluation
of the objective function and, as we propose in this paper, of its
gradient) is performed only at the point denoted by number 2.

Let us suppose that the interval shown in light grey in
Fig.~\ref{fig_OnePoint}a is chosen for the further partitioning.
Thus, at the third iteration, three smaller subintervals are
generated (see Fig.~\ref{fig_OnePoint}b). As one can see from
Fig.~\ref{fig_OnePoint}c, the trial point of the fourth iteration
coincides with the point 3 at which the trial has already been
executed. Therefore, there is no need to perform a new (costly)
evaluation of $f(x)$ and $f'(x)$ at this point, since the values
obtained at the previous iteration can be used. These values can
be stored in a specially designed vertex database and is simply
retrieved on demand without re-evaluations of the functions. For
example, Fig.~\ref{fig_OnePoint}d illustrates the situation after
12 iterations. It can be seen from this figure that $23$ intervals
have been generated by only 9 trial points.

Now we can describe the general scheme of a hyperinterval
partitioning by assuming (without loss of generality) that the
search hyperinterval $D$ in~(\ref{LGOP_D}) is an $N$-dimensional
hypercube and the first trial is performed at the vertex~$a$ (the
scheme starting from the vertex $b$ is obtained analogously). Let
a hyperinterval $D_t=[a_t, b_t]$ of a current partition $\{D^k\}$
of $D=[a,b]$ be chosen for partitioning at an iteration $k \geq 1$
of the algorithm. The operation of partitioning the selected
hyperinterval~$D_t$ is performed as follows.

\vspace{2mm}

\begin{description}
 \item[{\bf Step 1.}] Determine points $u$ and $v$ by the following formulae
 \be
  u=(a(1),\ldots,a(i-1),a(i) + \frac{2}{3}
  (b(i)-a(i)),a(i+1),\ldots,a(N)), \label{partition:u}
 \ee
 \be
  v=(b(1),\ldots,b({i-1}),b(i) + \frac{2}{3}
  (a(i)-b(i)),b(i+1),\ldots,b(N)),\label{partition:v}
 \ee
where $a(j)=a_t(j),\ b(j)=b_t(j),\ 1 \leq j \leq N$, and $i$ is
given by the equation
 \be
   i = \arg \min\, \max_{1 \leq j \leq N} | b(j) - a(j) |.
   \label{side_i}
 \ee
Get (evaluate or read from the vertex database) the values of the
objective function $f(x)$ and its gradient $f'(x)$ only at the
point $u$.

 \item[{\bf Step 2.}] Divide the hyperinterval $D_t$ into three
hyperintervals of equal volume by two parallel hyperplanes that
are perpendicular to the longest edge $i$ of $D_t$ and pass
through the points $u$ and $v$.

\hspace{2mm}The hyperinterval $D_t$ is so substituted by three new
hyperintervals with indices $t'=t$, $m +1$, and $m + 2$ (where $m
= m(k)$ is  the number of hyperintervals at the beginning of the
iteration $k$) determined by the vertices of their main diagonals
 \be \label{partition:D_t}
   a_{t'} = a_{m+2} = u, \  b_{t'} = b_{m+1} = v,
 \ee
 \be \label{partition:D_1}
   a_{m+1} = a_t, \ b_{m+1} = v,
 \ee
 \be \label{partition:D_2}
   a_{m+2} = u, \  b_{m+2} = b_t.
 \ee

Augment the current number of hyperintervals $m$ by 2.

\end{description}

From the partitioning scheme described above it can be observed
that, contrary to many traditional partitioning strategies (see,
e.g., \cite{Evtushenko&Posypkin(2011), Gergel(1997),
Jones:et:al.(1993), Kvasov:et:al.(2003), Pinter(1996)}), the
condition
 $$
  a_i(j) < b_i(j) \hspace{2mm} \forall j: j=1,\ldots,N,
 $$
does not have to be satisfied for all hyperintervals $D_i \subset
D$, and their main diagonals determined by the vertices $a_i$ and
$b_i$ can be oriented in different ways. However, as theoretically
shown in~\cite{Sergeyev&Kvasov(2008), Sergeyev(2000)}, the
hyperintervals orientations are not arbitrary and a special
linking of hyperintervals generated at different iterations can be
established with some efforts.

This smart linking will allow us to store information about
vertices and the corresponding values of $f(x)$ and $f'(x)$ in a
special database, thereby avoiding redundant functions
evaluations. The objective function and its gradient will be
calculated at a vertex only once, stored in the database, and read
when required. The new partitioning strategy generates trial
points in such a regular way that one vertex where the functions
are evaluated can belong to several (up to~$2^N$) hyperintervals
(see, for example, a trial point at the 8-th iteration in
Fig.~\ref{fig_OnePoint}d). Therefore, the time-consuming operation
of the functions evaluations is replaced by a significantly faster
operation of reading (up to $2^N$ times) the functions values from
the database. In this way, the new partitioning strategy
considerably speeds up the search, especially when problems of
high dimensions are considered (see~\cite{Kvasov&Sergeyev(2003),
Sergeyev(2000), Sergeyev(2005)}).

Note also that the possibility to choose the sequence of trial
points among either the points $a_i$ or the points $b_i$ (or among
other   $2^N-2$ vertices) of hyperintervals $D_i$ (see Step 1 of
the scheme) offers an important tool for accelerating the global
search when some additional information about the objective
function is known (we will see an example of this situation in
Section~4). Note the center-sampling partitioning strategies (see,
e.g.,~\cite{Evtushenko&Posypkin(2011), Gablonsky&Kelley(2001),
Jones:et:al.(1993)}) do not have this property.

\subsection{Lower bounding}

Let us consider an iteration $k \geq 1$ of the new algorithm and a
current partition $\{D^k\}$ of the search hyperinterval $D=[a,b]$
into hyperintervals $D_i=[a_i, b_i]$, $1 \leq i \leq m(k)$; over
these hyperintervals the values of both the function and its
gradient are obtained (evaluated or read from the vertex database)
at trial points $x^{j(k)}=a_i$, $j(k) \geq 1$. In order to choose
some hyperintervals for the further partition, the goodness
(expressed by the so-called \textit{characteristic}, see, e.g.,
\cite{Pinter(1996), Sergeyev&Kvasov(2008), Sergeyev&Kvasov(2011),
Strongin&Sergeyev(2000)}) of the hyperintervals with respect to
the global search is estimated by the algorithm. Better is the
characteristic of a hyperinterval (in some predetermined sense),
higher is the possibility to find the global minimizer within this
hyperinterval. This hyperinterval is, therefore, a good candidate
for a subdivision at the next iteration of the algorithm.

An estimate of the lower bound of $f(x)$ over a hyperinterval is
one of the possible characteristics of this hyperinterval. The
following result holds.

\begin{thm} \label{Theorem1}
Let $\tilde{K}$ be an estimate of the Lipschitz constant $K$ for
$f'(x)$ from (\ref{LGOP_K}), $\tilde{K} \ge K$ and $D_i = [a_i,
b_i]$ be a hyperinterval of a current partition $\{D^k\}$ with a
trial point $a_i$. Then, a value $R_i(\tilde{K})$ of the
characteristic of $D_i$ can be found such that it is the lower
bound of $f(x)$ over~$D_i$, i.e., $R_i(\tilde{K}) \leq f(x)$, $x
\in D_i$.
\end{thm}

\begin{pf} Let us prove the theorem in a constructive way.
It is known (see, e.g.,~\cite{Evtushenko&Posypkin(2011),
Nesterov(2004), Nocedal&Wright(1999)}) that for a differentiable
function $f(x)$ over a hyperinterval $D_i = [a_i, b_i]$ the
following inequality is satisfied:
 \be
  f(x) \geq Q(x,\tilde{K}), \hspace{3mm} x \in D_i,
  \label{Qmin1}
 \ee
where the quadratic minorant function $Q(x,\tilde{K})$ is defined
over $D_i$ as
 \be
  Q(x,\tilde{K}) = f(a_i) + \langle f'(a_i),(x-a_i)\rangle -
  0.5\tilde{K}\| x - a_i \|^2, \hspace{3mm} x \in
  D_i.
  \label{Qmin2}
 \ee
Here $\langle \cdot,\cdot\rangle$ is the scalar product, $\|\cdot
\|$ is the Euclidean norm in $R^N$, and
 $$
  g(x) = f(a_i) - \langle f'(a_i),(x-a_i)\rangle
 $$
is the linear approximation of $f(x)$ over $D_i$.

From inequality~(\ref{Qmin1}) the following estimates can be
obtained:
 $$
  f(x) \geq f(a_i) + \langle f'(a_i),(x-a_i)\rangle -
  0.5\tilde{K}\|b_i - a_i \|^2 \geq
 $$
 $$
  \geq F_i - 0.5\tilde{K}\|b_i - a_i \|^2, \hspace{3mm} x \in D_i,
 $$
where $F_i$ is the minimum value of the linear approximation
$g(x)$ over $D_i$, i.e.,
 \be \label{Fi1}
   F_i = f(a_i) + \min_{x \in D_i}\langle f'(a_i),(x-a_i)\rangle.
 \ee

Since the function $g(x)$ is linear, its minimum~(\ref{Fi1}) is
obtained in   the vertex $z_i$ of the hyperinterval $D_i=[a_i,
b_i]$ which coordinates $z_i(j), j=1, \ldots, N$, can be
calculated as follows:
 \be
   z_i(j) = \left\{
 \begin{array}{ll}
   a_i(j), & {\rm if \hspace{2mm} either \hspace{2mm} } b_i(j) > a_i(j)
                     {\rm\hspace{2mm} and \hspace{2mm}}
             \frac{\partial f(a_i)}{\partial x(j)} \geq 0, \\
           & {\rm or \hspace{2mm} } b_i(j) < a_i(j)
                     {\rm\hspace{2mm} and \hspace{2mm}}
             \frac{\partial f(a_i)}{\partial x(j)} < 0; \\

   b_i(j), & {\rm if \hspace{2mm} either \hspace{2mm} } b_i(j) > a_i(j)
                     {\rm\hspace{2mm} and \hspace{2mm}}
             \frac{\partial f(a_i)}{\partial x(j)} < 0, \\
           & {\rm or \hspace{2mm}} b_i(j) < a_i(j)
                     {\rm\hspace{2mm} and \hspace{2mm}}
             \frac{\partial f(a_i)}{\partial x(j)} \geq 0. \\
 \end{array}
 \right.
 \label{zj}
 \ee
The corresponding value $F_i$ from~(\ref{Fi1}) is therefore equal
to
 \be \label{Fi2}
   F_i = f(a_i) + \langle f'(a_i),(z_i-a_i)\rangle.
 \ee
It is clear now that the value
 \be
   R_i = R_i(\tilde{K}) = F_i - 0.5\tilde{K}\|b_i-a_i\|^2
   \label{Ri}
 \ee
satisfies the inequality
 $$
   R_i \leq f(x), \hspace{3mm} x \in D_i,
 $$
and, therefore, it can be taken as the characteristic value of
$D_i$ that estimates the lower bound of $f(x)$ over $D_i$. The
theorem has been proved. \hfill\rule{5pt}{5pt}
\end{pf}

Note that analogous results can be obtained in the case of
hyperintervals~$D_i$ with trial points $b_i$ rather than $a_i$.

\begin{figure}[t]
\centerline{\psfig{file=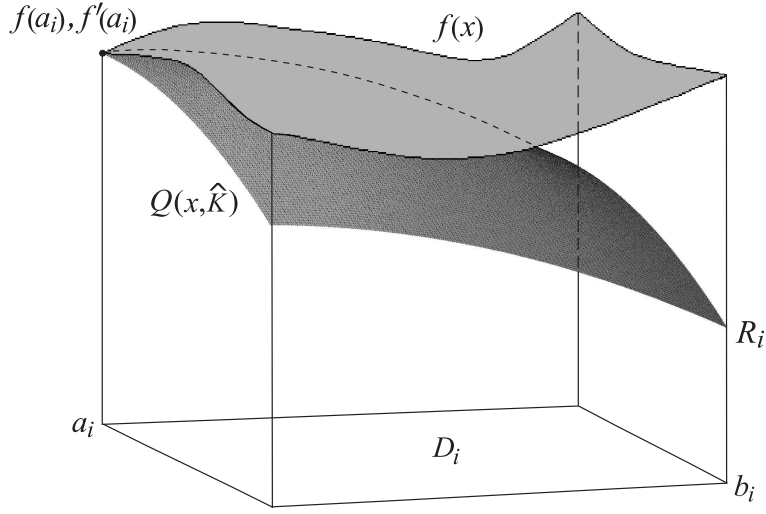,width=100mm,height=68mm,angle=0,silent=}}
\vspace*{2mm}\caption{A quadratic minorant function
$Q(x,\tilde{K})$ for $f(x)$ over a hyperinterval $D_i=[a_i, b_i]$}
\label{Fig:Ri}
\end{figure}

In Fig.~\ref{Fig:Ri}, a quadratic minorant function
$Q(x,\tilde{K})$ from~(\ref{Qmin2}) is illustrated for $f(x)$ over
a hyperinterval $D_i$. Here, the characteristic value $R_i$
coincides with the minimum value of $Q(x,\tilde{K})$ obtained at
the point $b_i$ of the main diagonal of $D_i$. In general, as it
can be seen from~(\ref{Qmin2}), the value $R_i$ is smaller than or
equal to the minimum value of $Q(x,\tilde{K})$ over $D_i$.

\subsection{Nondominated hyperintervals and their graphical representation}

By using the obtained characteristics of hyperintervals, the
relation of domination can be established between every two
hyperintervals of a current partition $\{D^k\}$ of $D$ and a set
of nondominated hyperintervals can be identified for a possible
subdivision at the current iteration of the new algorithm
(see~\cite{Kvasov&Sergeyev(2009), Sergeyev&Kvasov(2006)}).

 \begin{dfn} \label{Def:Dominated}
Given an estimate $\tilde{K}>0$ of the Lipschitz constant~$K$
from~(\ref{LGOP_K}), a hyperinterval $D_i=[a_i,b_i]$
\textit{dominates} a hyperinterval $D_j=[a_j, b_j]$ with respect
to $\tilde{K}$ if
 $$
  R_i(\tilde{K}) < R_j(\tilde{K}).
 $$
  \end{dfn}
  \begin{dfn} \label{Def:NonDominated}
A hyperinterval $D_t=[a_t, b_t]$ is said to be
\textit{nondominated with respect to} $\tilde{K}>0$ if for the
chosen value $\tilde{K}$ there is no other hyperinterval
in~$\{D^k\}$ which dominates $D_t$.
 \end{dfn}

\begin{figure}[!t]
\centerline{\epsfig{file=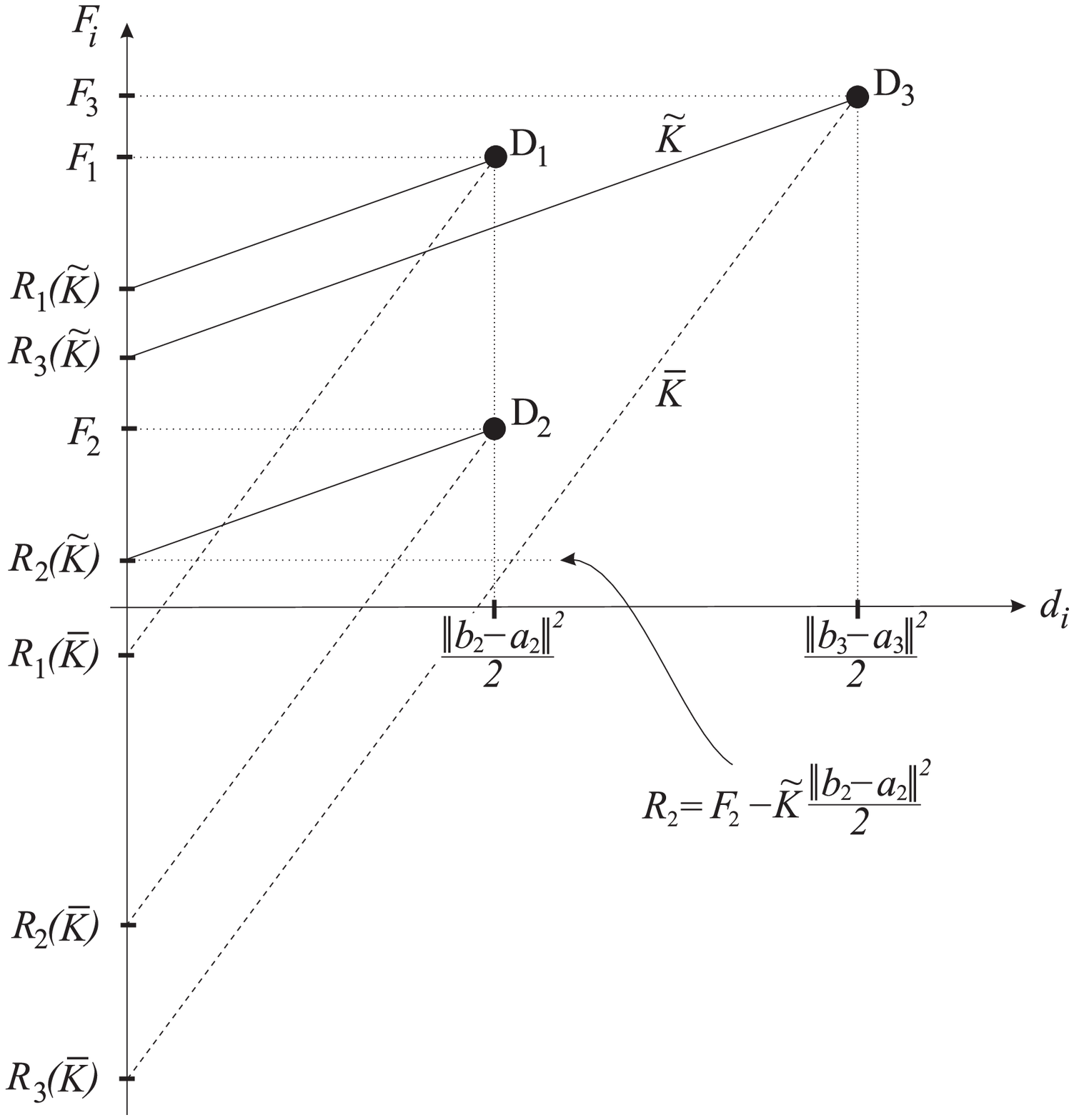,width=100mm,height=102mm,angle=0,silent=}}
\caption{Graphical representation of hyperintervals}
\label{Fig:Diagram}
\end{figure}

Let us now generalize the approach proposed by the authors
in~\cite{Kvasov&Sergeyev(2009)} for the one-dimensional prototype
and show that both a multi-dimensional interval $D_i=[a_i, b_i]$
of a current partition $\{D^k\}$ and the respective characteristic
$R_i$ using the gradient can be represented in a two-dimensional
diagram similar to those proposed in~\cite{Jones:et:al.(1993),
Sergeyev&Kvasov(2006)} for derivative free methods. Difficulties
in the construction of such a diagram were among the main reasons
that prevented people to propose methods using several estimates
of $K$ in their work.

So, we take for the dot, corresponding to $D_i$, the vertical
coordinate $F_i$ from~(\ref{Fi1})--(\ref{Fi2}) and the horizontal
coordinate $d_i$ equal to half of the squared length of the main
diagonal of $D_i$, i.e.,
 $$
   d_i = 0.5 \|b_i - a_i\|^2.
 $$

For example, in Fig.~\ref{Fig:Diagram}, a partition of the search
domain $D$ consisting of three hyperintervals is represented by
the dots $D_1$, $D_2$, and $D_3$. Let us suppose that the
Lipschitz constant~$K$ for the gradient $f'(x)$ is estimated
by~$\tilde{K}$, $\tilde{K} \geq K$. The characteristic $R_i$ of a
hyperinterval $D_i$, $i=1,2,3$, can be graphically obtained as the
vertical coordinate of the intersection point of the line passed
through the point $D_i$ with the slope $\tilde{K}$ and the
vertical coordinate axis (see Fig.~\ref{Fig:Diagram}). It is easy
to see, that with respect to the estimate $\tilde{K}$ the
hyperinterval $D_2$ dominates both hyperintervals $D_1$ and $D_3$
and the hyperinterval $D_3$ dominates $D_1$.

If a higher estimate $\bar{K} > \tilde{K}$ of the Lipschitz
constant $K$ is considered (see Fig.~\ref{Fig:Diagram}), the
hyperinterval $D_2$ still dominates $D_1$ with respect to
$\bar{K}$, because $R_2 (\bar{K}) < R_1(\bar{K})$. But $D_2$ in
its turn is dominated by the hyperinterval~$D_3$ with respect to
$\bar{K}$, because $R_2(\bar{K})
> R_3 (\bar{K})$ (see Fig.~\ref{Fig:Diagram}).

Since the exact Lipschitz constant $K$ for $f'(x)$ (or its valid
overestimate) is unknown in the stated problem, the following
definition can be useful.

\begin{dfn} \label{DefNonDomination}
A hyperinterval $D_t \in \{D^k\}$ is called \textit{nondominated}
if there exists an estimate $0 < \tilde{K} <\infty$ of the
Lipschitz constant~$K$ such that $D_t$ is nondominated with
respect to $\tilde{K}$.
\end{dfn}

This means that nondominated hyperintervals are those with the
smallest characteristics~(\ref{Ri}) for some particular estimate
of the Lipschitz constant for the gradient $f'(x)$. For example,
in Fig.~\ref{Fig:Diagram} the hyperintervals $D_2$ and $D_3$ are
nondominated.

\begin{figure}
\centerline{\psfig{file=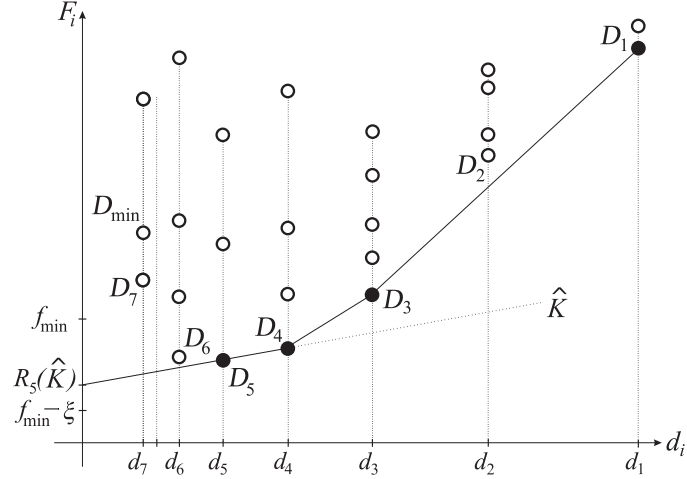,width=0.70\textwidth,height=0.50\textwidth,angle=0,silent=}}
\caption{The two-dimensional diagram representing dominated (white
dots) and nondominated (black dots) hyperintervals of a current
partition of the search domain} \label{Fig:Optimal}
\end{figure}

It can be demonstrated following the reasoning used in
~\cite{Kvasov&Sergeyev(2009), Sergeyev&Kvasov(2006)} that
nondominated hyperintervals (in the sense of
Def.~\ref{DefNonDomination}) are located on the lower-right convex
hull of the set of dots representing the hyperintervals of the
current partition of $D$ and can be efficiently found by applying
  algorithm for identifying the convex hull of the dots (see,
e.g.,~\cite{He:et:al.(2002), Jones:et:al.(1993),
Sergeyev&Kvasov(2008)}). In Fig.~\ref{Fig:Optimal}, the
hyperintervals represented by the dots $D_1$ (the largest
hyperinterval), $D_3$, $D_4$, and $D_5$ are nondominated
hyperintervals.

It has been shown in~\cite{Sergeyev&Kvasov(2006)} that the
hyperintervals of a current partition of~$D$ form several groups
characterized by the length of their main diagonals. The
hyperintervals from a group are represented graphically by dots
with the same horizontal coordinate. For example, in
Fig.~\ref{Fig:Optimal} there are seven different groups of
hyperintervals with the horizontal coordinates equal to $d_i$,
$i=1, \ldots, 7$, and one empty group (with the horizontal
coordinate between $d_7$ and $d_6$ in Fig.~\ref{Fig:Optimal}).
Empty groups correspond to hyperintervals which are not present in
the current partition but can be generated (or were generated) at
the successive (previous) iterations of the method.

As demonstrated in~\cite{Sergeyev&Kvasov(2006), Sergeyev(2000)}, a
correspondence between the length of the main diagonal of a
hyperinterval $D_i$ and a non-negative integer number can be
established, which indicates the number of subdivisions of the
initial domain~$D$ necessary to obtain the hyperinterval~$D_i$. At
each iteration $k \geq 1$ this number can be considered as an
index $s(k)$ of a group of equal hyperintervals where
 \be \label{GroupIndex}
  0 \leq q_\infty(k) \leq s(k) \leq q_0(k) < +\infty
 \ee
and $q_\infty (k)$ and $q_0(k)$ are indices corresponding to the
groups of the largest and smallest hyperintervals of the current
partition of $D$, respectively (for example, in
Fig.~\ref{Fig:Optimal}, $q_0(k)= q_\infty(k)+7$). During
partitioning, diagonals of hyperintervals become smaller, while
the corresponding group indices grow up consecutively starting
from $q_\infty(1)=0$ (see~\cite{Sergeyev&Kvasov(2006)} for
details).

Once a nondominated hyperinterval $D_t=[a_t, b_t]$ is determined
(with respect to some estimate $\tilde{K}$ of the Lipschitz
constant $K$), it can be subdivided at the next iteration of the
algorithm if the following condition is satisfied:
 \be  \label{ConditionOfImprovement}
   R_t(\tilde{K}) \leq f_{\rm min} (k) - \xi,
 \ee
where $R_t$ is calculated by~(\ref{Ri}), $f_{\rm min} (k)$ is the
record value, i.e., the current minimal function value (attained
at the record point $x_{\rm min} (k)$), and $\xi$ is the parameter
of the algorithm, $\xi \geq 0$ (it can be set in different ways,
see Section~\ref{sectionResults}). Notice that both the record
value and the record point can be changed after performing some
better trial during partitioning, but the record value remains
always greater than or equal to the vertical coordinate of the
lowest dot (dot~$D_5$ in Fig.~\ref{Fig:Optimal}).

Condition~(\ref{ConditionOfImprovement}) prevents the algorithm
from subdividing already well-explored small hyperintervals. For
example, among nondominated hyperintervals in
Fig.~\ref{Fig:Optimal} (black dots), the hyperinterval $D_5$ does
not satisfy this condition and therefore is excluded from being
partitioned at the next iteration of the method.

It should be mentioned in this occasion that, together with
nondominated hyperintervals, a hyperinterval $D_{\rm min}(k) =
[a_{\rm min}, b_{\rm min}]$ containing the record point (called
hereafter the record hyperinterval) is also considered for a
possible partition during the work of the algorithm as it will be
explained in the next Section. Among different hyperintervals
 the record point~$x_{\rm min}(k)$ can belong to (up to
$2^N$), the record hyperinterval is that with the smallest
characteristic and can be changed during subdivisions. In
Fig.~\ref{Fig:Optimal}, the record hyperinterval is represented by
the dot $D_{\rm min}$ (note that this dot can be not the lowest
one, as in Fig.~\ref{Fig:Optimal}). Hereafter, the index of the
group the hyperinterval $D_{\rm min}(k)$ belongs to will be
indicated as $p(k)$ (during the work of the algorithm the
satisfaction of inequalities~(\ref{GroupIndex}) is ensured for
this index which can be eventually updated together with $q_0(k)$
and $q_\infty(k)$; see~\cite{Sergeyev&Kvasov(2006)} for details).
In Fig.~\ref{Fig:Optimal}, $p(k) = q_0(k)$ and, therefore, the
hyperinterval $D_{\rm min}(k)$ is among the smallest
hyperintervals of the illustrated partition of $D$.

\section{New Algorithm}

In this Section, the new algorithm for solving problem
(\ref{LGOP_f})--(\ref{LGOP_D}) is described. First, the new method
is presented and its computational scheme is given, then its
convergence properties are analyzed.

The new algorithm consists of the following explicitly defined
phases: (1) an exploration phase, at which an examination of large
hyperintervals (possibly located far away from the record point)
is performed in order to capture new subregions with better
function values; (2) a record improvement phase, at which the
algorithm tries to better inspect the subregion around the record
point. Several subdivisions of different hyperintervals can be
performed at a single iteration of the new method (this more
general notion of an iteration with respect to that of
Section~\ref{sectionNewStrategy} is often used in the Lipschitz
global optimization algorithms with multiple estimates of
Lipschitz constants, see, e.g.,~\cite{Gablonsky&Kelley(2001),
Jones:et:al.(1993), Kvasov&Sergeyev(2009),
Sergeyev&Kvasov(2006)}).

The exploration phase consists of several iterations (namely,
$N+1$ where $N$ is the problem dimension), each serves for
determining nondominated hyperintervals and partitioning them.
Since each subdivision of a hyperinterval by the
scheme~(\ref{partition:u})--(\ref{partition:D_2}) is performed
perpendicularly to only one side of the hyperinterval (to the
longest side from~(\ref{side_i})), the number of iterations within
a phase of the algorithm should be correlated with the
hyperintervals dimension.

This phase is interrupted after finishing an iteration if an
improvement on at least $1\%$ of the minimal function value is
reached, i.e., if
   \be \label{ConditionOfLocalImprovement}
     f_{\rm min}(k) \leq f_{\rm min}^{\rm prec}- 0.01 | f_{\rm min}^{\rm prec}
     |,
   \ee
where $f_{\rm min}^{\rm prec}$ is the record value memorized at
the start of the exploration phase.

Condition~(\ref{ConditionOfLocalImprovement}) is verified after
each iteration of the exploration phase and is used to switch the
algorithm to the record improvement phase. This local phase is
also launched when the exploration phase finishes without having
improved the record value, but only if the record hyperinterval
$D_{\rm min} (k)$ is not the smallest one within the current
partition of hyperintervals (for example, in
Fig.~\ref{Fig:Optimal}, the record hyperinterval is among the
smallest hyperintervals). Otherwise, the algorithm re-initiates
another global exploration phase without forcing the local one.

The record improvement phase reflects the already well-established
fact in global optimization affirming the benefits of the record
improvement during the global search (see, e.g., the references
given in~\cite{Horst&Pardalos(1995), Kvasov&Sergeyev(2009),
Lera&Sergeyev(2010a), Sergeyev&Kvasov(2008)}). At a single
iteration, it performs several subdivisions (namely, $N$) of the
record hyperinterval trying to improve the record value. During
this process a new record value can appear. In this case, a new
record hyperinterval can be considered for remaining subdivisions.

The record hyperinterval subdivisions are performed by means of
the   one-point-based strategy described in
Section~\ref{sectionNewStrategy}. Of course, other possible local
improvement techniques can be used for this scope (see,
e.g.,~\cite{Gaviano&Lera(2002), Nesterov(2004),
Nocedal&Wright(1999)}) but in this case the resulting trial points
can not be managed within the vertex database mentioned in
Section~\ref{sectionNewStrategy}.

It is important that the available gradient information allows us
to terminate automatically the record improvement phase. In fact,
the record hyperinterval is not further subdivided when the
gradient projection on the directions parallel to the record
hyperinterval sides becomes non-negative, i.e., when the following
condition is satisfied:
 \be
   \frac{\partial f(a_{\rm min})}{\partial x(j)}(b_{\rm min}(j)-a_{\rm min}(j)) \geq 0
   \hspace{3mm}\forall j: j=1,\ldots,N.
 \label{ExtremumCond}
 \ee

Either in this case or when the prefixed number $N$ of
subdivisions are normally performed (without meeting
conditions~(\ref{ExtremumCond})), the algorithm is switched again
to the global exploration phase and continues its work.

The algorithm stops when the number of generated trial points
reaches the maximal allowed number $P_{\rm max}$. The satisfaction
of this termination criterion is verified after every subdivision
of a hyperinterval. The current record value~$f_{\rm min}$ and the
current record point $x_{\rm min}$ can be taken as approximations
of the global minimum value $f^*$ and the global minimizer $x^*$
from~(\ref{LGOP_f}), respectively.

\vspace*{3mm}

A formal description of the new algorithm   follows below (we
assume without loss of generality that the admissible region
$D=[a,b]$ in~(\ref{LGOP_D}) is an $N$-dimensional hypercube).

\begin{description}
\item {\bf Step 0} (\textit{Initialization}). Set the iteration
counter $k:=1$. Let the first evaluation of $f(x)$ and $f'(x)$ be
performed at the vertex $a$ of the initial hyperinterval
$D=[a,b]$, i.e., $x^1 :=a$. Set the current partition of the
search interval as $D^1 := \{[a_1,b_1]\}$, where $a_1=a$, $b_1=b$,
and the current number of hyperintervals $m(1):=1$. Set $f_{\rm
min} (1) := f(x^1)$, $x_{\rm min}(1):=a$, and $D_{\rm
min}(1):=[a_1,b_1]$. Set group indices
$q_\infty(1):=q_0(1):=p(1):=0$.

\end{description}

Suppose now that $k \geq 1$ iterations of the algorithm have
already been executed. The next iterations of the algorithm
consist of the following steps.

 \begin{description}

\item {\bf Step 1} (\textit{Exploration Phase}). Memorize the
current record $f_{\rm min}^{\rm prec}:=f_{\rm min}(k)$, set the
counter of iterations during the exploration phase $k_g:=1$ and
perform the following steps:

 \begin{description}

  \item[{\bf Step 1.1.}] Identify the set of nondominated
hyperintervals considering only groups of large hyperintervals
(namely, those with the current indices from $q_\infty(k)$ up to
$\lceil(q_\infty(k)+p(k))/2\rceil$). Subdivide those nondominated
hyperintervals which satisfy
inequality~(\ref{ConditionOfImprovement}) and produce new trial
points (or read the existing ones from the vertex database)
according to Section~\ref{sectionNewStrategy}. Set $k:=k+1$ and
update hyperintervals indices if necessary
(see~\cite{Sergeyev&Kvasov(2006)} for details).

 \item[{\bf Step 1.2.}] If
condition~(\ref{ConditionOfLocalImprovement}) is satisfied, then
go to Step 2 and execute the record improvement phase. Otherwise,
go to Step 1.3.

 \item[{\bf Step 1.3.}] Increase the counter $k_g:= k_g +1$: check
whether $k_g \leq N$. If this is the case, then go to Step 1.1
(continue the exploration of large hyperintervals). Otherwise, go
to Step 1.4 (perform the final iteration of the exploration phase
by considering more hyperintervals groups).

 \item[{\bf Step 1.4.}] Identify the set of nondominated
hyperintervals considering the current groups of hyperintervals
from $q_\infty(k)$ up to $p(k)$. Subdivide those nondominated
hyperintervals which satisfy
inequality~(\ref{ConditionOfImprovement}) and produce new trial
points (or read the existing ones from the vertex database)
according to Section~\ref{sectionNewStrategy}. Set $k:=k+1$,
update all necessary indices.

 \item[{\bf Step 1.5.}] If the record hyperinterval is not the
smallest one, i.e., if $p(k) < q_0(k)$, then then go to Step 2 and
execute the record improvement phase. Otherwise, go to Step 1 and
repeat the exploration phase updating the value $f_{\rm min}^{\rm
prec}$.

 \end{description}

\item {\bf Step 2} (\textit{Record Improvement Phase}). Set
$k:=k+1$. Set the counter of iterations during the record
improvement phase $k_l:=1$ and perform the following steps:

 \begin{description}

 \item[{\bf Step 2.1.}] Subdivide the record hyperinterval $D_{\rm
 min}(k)$ and produce a new trial
point (or read the existing one from the vertex database)
according to Section~\ref{sectionNewStrategy}. Update
hyperintervals indices and the record hyperinterval index if
necessary.

 \item[{\bf Step 2.2.}] Increase the counter $k_l:= k_l +1$: check
whether $k_l \leq N$. If this is the case, then go to Step 1
(perform a new exploration of large hyperintervals). Otherwise, go
to Step 2.1 (continue the local exploration of the subregion near
to the record point).

 \end{description}

\end{description}

Let us now study convergence properties of the new method during
minimization of the function $f(x)$
from~(\ref{LGOP_f})--(\ref{LGOP_D}) when the maximal allowed
number of generated trial points $P_{\rm max}$ is equal to
infinity. In this case, the algorithm does not stop (the number of
iterations $k$ goes to infinity) and an infinite sequence of trial
points~$\{x^{j(k)}\}$ is generated.

\begin{dfn}
The convergence of an infinite sequence of trial points
$\{x^{j(k)}\}$ generated by a global optimization method is called
\textit{everywhere dense} if for any point $x \in D$ and any
$\delta >0$ there exist an iteration number $k(\delta) \geq 1$ and
a point $x' \in \{x^{j(k)}\}$, $k > k(\delta)$, such that $\|x -
x'\| < \delta$.
\end{dfn}

\begin{thm} \label{TheoremConvergence}
The new algorithm manifests the everywhere dense convergence.
\end{thm}

\begin{pf} Every subdivision of a
hyperinterval by the partitioning scheme from
Section~\ref{sectionNewStrategy} produces three new hyperintervals
with the same volume equal to the third part of the volume of the
subdivided hyperinterval and smaller main diagonals. Trial points
generated by the new algorithm are at one of the vertices of each
generated hyperinterval. Therefore, fixed a positive value
of~$\delta$, it is sufficient to demonstrate that after a finite
number of iterations~$k(\delta)$, the largest hyperinterval of the
current partition of the search domain $D$ will have the length of
its main diagonal smaller than $\delta$. In this case, in
$\delta$-neighborhood of any point of $D$ there will exist at
least one trial point generated by the algorithm.

Let us fix an iteration $k'$ of the method and consider the group
$q_\infty(k')$ of the largest hyperintervals of the partition
$\{D^{k'}\}$ on its two-dimensional graphic representation. This
group is always taken into account when nondominated
hyperintervals are looked for at the exploration phase of the
algorithm. As it follows from~Def.~\ref{DefNonDomination}, a
hyperinterval $D_t \in \{D^{k'}\}$ from this group with the
corresponding smallest value $F_t$ from~(\ref{Fi1})--(\ref{Fi2}),
must be partitioned and substituted by three smaller
hyperintervals at the current iteration of the algorithm because
it is a nondominated hyperinterval and
condition~(\ref{ConditionOfImprovement}) is satisfied for it.

Since each group consists of a finite number of hyperintervals,
after a sufficiently large number of iterations $k
> k'$ all hyperintervals of the group~$q_\infty(k)$ of the largest hyperintervals
will be subdivided. The group index~$q_\infty(k)$ will increase
and the same procedure will be repeated with a new group of the
largest hyperintervals, thus making the largest hyperintervals
smaller and smaller.

It can be also noted that the record hyperinterval $D_{\rm min}$
is itself represented by a dot in the two-dimensional diagram of
the current partition. It can be subdivided either separately
during the record improvement phase, or as a nondominated
hyperinterval during the exploration phase at which the
satisfaction of condition~(\ref{ExtremumCond}) is not taken in
consideration.

Thus, there exists a finite number $k(\delta)$ such that after
executing~$k(\delta)$ iterations of the algorithm the largest
hyperinterval of the current partition $\{D^{k(\delta)}\}$ will
have the main diagonal smaller than $\delta$.
\hfill\rule{5pt}{5pt}
\end{pf}

To conclude the theoretical study of the new algorithm we would
like to highlight that the usage of all possible estimates of the
Lipschitz constant in its work leads to the convergence of the
everywhere dense type.  If the Lipschitz constant $L$ (or its
valid estimate) of the objective function $f(x)$ or the Lipschitz
constant $K$ (or its valid estimate) of the gradient $f'(x)$ can
be used by a global optimization method, other types of
convergence can be established for such an algorithm (see, e.g.,
methods from~\cite{Horst&Tuy(1996), Pinter(1996),
Sergeyev&Kvasov(2008), Sergeyev&Kvasov(2011),
Strongin&Sergeyev(2000)}).

\section{Numerical results} \label{sectionResults}

In this Section, we present numerical results performed to compare
the new algorithm with two methods belonging to the same class of
the one-point-base partitioning methods: the DIRECT algorithm
from~\cite{Jones:et:al.(1993)} and its locally-biased modification
DIRECT{\it l} \hspace{1mm}from \cite{Gablonsky&Kelley(2001)}. Both
of them use the center-sampling partitioning strategy and work
with a set of Lipschitz constants for the objective function
$f(x)$ from~(\ref{LGOP_f}). The implementation of these two
methods (downloadable from
\url{http://www4.ncsu.edu/~ctk/SOFTWARE/DIRECTv204.tar.gz}) has
been used in all the experiments following the way of the
multicriteria  comparison proposed
in~\cite{Sergeyev&Kvasov(2006)}.

In order to make easier the numerical comparison with the
DIRECT-based algorithms, the value~$\xi$
from~(\ref{ConditionOfImprovement}) was set as in the DIRECT
method,~i.e.,
 \be  \label{DirEpsilon}
    \xi = \epsilon | f_{\rm min} (k) |, \hspace{3mm} \epsilon
    \geq 0.
 \ee
The recommended value of $\epsilon=10^{-4}$
(see~\cite{Jones:et:al.(1993), Sergeyev&Kvasov(2006)}) was used
in~(\ref{DirEpsilon}).

In accordance with~\cite{Sergeyev&Kvasov(2006)}, the global
minimizer $x^* \in D$ was considered to be found when a method
generated a trial point~$x'$ inside a hyperinterval with a vertex
$x^*$ and the volume smaller than the volume of the initial
hyperinterval $D=[a,b]$ multiplied by an accuracy
coefficient~$\Delta$, $0< \Delta \leq 1$, i.e.,
 \be \label{x*Found}
  |x'(i) - x^*(i) | \leq \sqrt[N]{\Delta}(b(i)-a(i)), \hspace{3mm}
  1 \leq i \leq N,
 \ee
where $N$ is from~(\ref{LGOP_D}). The algorithm was stopped either
when the maximal number of trials $P_{\rm max}$ equal
to~1\,000\,000 was reached, or when condition~(\ref{x*Found}) was
satisfied (see~\cite{Sergeyev&Kvasov(2006)} for a discussion about
different stopping criteria in global optimization methods).

In our numerical experiments we used the same test classes, each
of 100 continuously differentiable functions, produced by the
GKLS-ge\-ne\-ra\-tor (see \cite{Gaviano:et:al.(2003)}) as
in~\cite{Sergeyev&Kvasov(2006)}. Particularly, eight GKLS D-type
classes of dimensions $N=2$, 3, 4, and 5 have been considered. For
each particular problem dimension~$N$ a `simple' and a `hard'
classes have been taken for the comparison
(see~\cite{Sergeyev&Kvasov(2006)} for a detailed description of
the classes).

For the convenience of the reader, we report here the four
criteria introduced in~\cite{Sergeyev&Kvasov(2006),
Sergeyev&Kvasov(2008)} that were used to compare the methods. The
following designations are required:

$P_s$ -- the number of trials performed by the method under
consideration to solve the problem number $s$, $1 \leq s \leq
100$, of a fixed test class.

$m_s$ -- the number of hyperintervals generated to solve the
problem $s$.

\vspace{1mm}

{\bf Criterion C1.} Number of trials $P_{s^*}$ required for a
method to satisfy condition~(\ref{x*Found}) for {\it all} 100
functions of a particular test class, i.e.,
 \be
   P_{s^*} = \max_{1 \leq s \leq 100} P_s, \hspace{5mm}
   s^* = \arg\max_{1 \leq s \leq 100} P_s. \label{f_num}
 \ee

{\bf Criterion C2.} The corresponding number of hyperintervals,
$m_{s^*}$, generated by the method, where $s^*$ is
from~(\ref{f_num}).

{\bf Criterion C3.} Average number of trials $P_{avg}$ performed
by the method during minimization of {\it all} 100 functions from
a particular test class, i.e.,
 \be
  P_{avg} = \frac{1}{100}\sum_{s=1}^{100} P_s. \label{avg_trials}
 \ee

{\bf Criterion C4.} Number $p$ (number $q$) of functions from a
class for which DIRECT or DIRECT{\it l} executed less (more)
function evaluations than the new algorithm. If~$P_s$ is the
number of trials performed by the new algorithm and $P_s'$ is the
corresponding number of trials performed by a competing method,
$p$ and $q$ are evaluated as follows
 \be
  p = \sum_{s=1}^{100} \sigma_s', \hspace{5mm} \label{winDIRECT}
  \sigma_s' =
  \left\{
    \begin{array}{ll}
       1, & P_s' < P_s,  \\
       0, & {\rm otherwise}.
   \end{array} \right.
 \ee
 \be
  q= \sum_{s=1}^{100} \sigma_s, \hspace{5mm} \label{winNEW}
  \sigma_s =
  \left\{
    \begin{array}{ll}
       1, & P_s < P_s',  \\
       0, & {\rm otherwise}.
   \end{array} \right.
 \ee

\vspace{1mm}

Results based on Criteria C1 and C2 are mainly influenced by
minimization of the most difficult functions of a class. Criteria
C3 and C4 deal with average data of a class. The number of
generated hyperintervals (Criterion C2) provides an important
characteristic of any partition algorithm for solving the
problem~(\ref{LGOP_f})--(\ref{LGOP_D}). In some way, it
corresponds to the qualitative examination of the search domain
$D$ during the work of the method. The greater is this number, the
more information about the behavior of the objective function is
available and, therefore, the smaller is the risk to miss its
global minimizer. Of course, algorithms should not generate many
redundant hyperintervals since this slows down the search and is
therefore a disadvantage of the method
(see~\cite{Sergeyev&Kvasov(2006)} for more details).

Results of numerical comparison of the methods with respect to
Criteria C1 and C2 with eight GKLS test classes are shown in
Tables \ref{table1}--\ref{table3}. The accuracy coefficient
$\Delta$ from~(\ref{x*Found}) is given in the second column of the
tables. Table~\ref{table1} reports the maximal number of trials
required for satisfying condition~(\ref{x*Found}) for half of the
functions of a particular class (columns ``50\%'') and for all 100
function of the class (columns ``100\%''). The notation `$>$
1\,000\,000 $(j)$' in Tables~\ref{table1} and \ref{table3} means
that after 1\,000\,000 function evaluations the method under
consideration was not able to solve~$j$ problems.
Table~\ref{table2} represents the ratio between the maximal number
of trials performed by DIRECT and DIRECT{\it l} with respect to
the corresponding number of trials performed by the new algorithm.
The numbers of generated hyperintervals (Criterion~C2) are
indicated in Table~\ref{table3}.

According to Tables \ref{table1} and \ref{table3}, the new
multidimensional algorithm requires much fewer trials than the
other two methods to ensure a thorough examination of the search
domain. Moreover, the advantage of the new method becomes even
more pronounced as the problem dimension grows or the problem
complexity increases.

In fact, on half of the test functions from each class (which were
the most simple for each method with respect to the other
functions of the class) the new algorithm already manifested a
very good performance with respect to DIRECT and DIRECT{\it l} in
terms of the number of generated trial points (see columns
``50\%'' in Table~\ref{table1}). When all the functions were taken
in consideration (and, consequently, difficult functions of the
class were considered too), the number of trials produced by the
new algorithm was much fewer in comparison with two other methods
(see columns ``100\%'' in Table~\ref{table1}), ensuring at the
same time a substantial examination of the admissible domain (see
Table~\ref{table3}).

\begin{table}[t]
\caption{Number of trial points for GKLS test functions (Criterion
C1).}
\begin{center} \scriptsize \label{table1}
\begin{tabular}{@{\extracolsep{\fill}}|c|c|c|r|r|r|r|r|r|}\hline
$N$ & $\Delta$ & Class &
\multicolumn{3}{c|}{50\%} & \multicolumn{3}{c|}{100\%}\\
\cline{4-9}& & &DIRECT &DIRECT{\it l} &New &DIRECT &DIRECT{\it l} &New\\
\hline
2 &$10^{-4}$ & simple &111  &152  &59   &1159  &2318 & 335\\
2 &$10^{-4}$ & hard   &1062 &1328 &182  &3201  &3414 & 1075\\
\hline
3 &$10^{-6}$ & simple &386   &591 &362 &12507 &13309 &2043\\
3 &$10^{-6}$ & hard   &1749  &1967 &416 &$>$1000000 (4) &29233 &2352 \\
\hline
4 &$10^{-6}$ & simple &4805  &7194  &2574 &$>$1000000 (4) &118744 &16976\\
4 &$10^{-6}$ & hard   &16114 &33147 &3773 &$>$1000000 (7) &287857 &20866\\
\hline
5 &$10^{-7}$ & simple &1660   &9246   &1757  &$>$1000000 (1)  &178217         &16300\\
5 &$10^{-7}$ & hard   &55092  &126304 &13662 &$>$1000000 (16) &$>$1000000 (4) &88459\\
\hline
\end{tabular}
\end{center}
\end{table}

\begin{table}[!t]
\caption{Improvement obtained by the new algorithm in terms of
Criterion C1.}
\begin{center}\scriptsize \label{table2}
\begin{tabular}{@{\extracolsep{\fill}}|c|c|c|c|c|}\hline
$N$ & $\Delta$ & {Class} &DIRECT/New &DIRECT{\it l}/New \\
\hline
2 &$10^{-4}$ & simple  & 3.46  & 6.92 \\
2 &$10^{-4}$ & hard  & 2.98  & 3.18 \\
\hline
3 &$10^{-6}$ & simple  & 6.12      & 6.51  \\
3 &$10^{-6}$ & hard  &$>$425.17  & 12.43 \\
\hline
4 &$10^{-6}$ & simple  &$>$58.91   & 6.99  \\
4 &$10^{-6}$ & hard  &$>$47.92   & 13.80 \\
\hline
5 &$10^{-7}$ & simple  &$>$61.35   & 10.93   \\
5 &$10^{-7}$ & hard  &$>$11.30   &$>$11.30 \\
\hline
\end{tabular}
\end{center}
\end{table}

Note also that maximal number of trials equal to 88459 (see
Table~\ref{table1}) required by the new method to solve all
problems of the hard five-dimensional class is obtained on the
function 5 of this class. If we use the new method with the
one-point-based strategy starting from the point $b$ rather than
from the point $a$ (see Section~\ref{sectionNewStrategy}), the
number of trials required by the new algorithm to solve this
particular problem becomes equal to $15238$. Thus, some a priori
knowledge on the objective function behavior can allow  us to
better select the vertex of the initial hyperinterval $D$ in which
the first trial will be executed  and, therefore, to accelerate
the search even more.

Table~\ref{table4} reports the average number of trials performed
during minimization of all 100 functions from the same GKLS
classes (Criterion~C3). The ``Improvement'' columns in these
tables represent the ratios between the average numbers of trials
performed by DIRECT and DIRECT{\it l} with respect to the
corresponding numbers of trials performed by the new algorithm.
The symbol~`$>$' reflects the situation when not all functions of
a class were successfully minimized by the method under
consideration in the sense of condition~(\ref{x*Found}). This
means that the method stopped when $P_{\rm max}$ trials had been
executed during minimization of several functions of this
particular test class. In these cases, the value of $P_{\rm max}$
equal to 1\,000\,000 was used in calculations of the average value
in~(\ref{avg_trials}), providing in such a way a lower estimate of
the average. As can be seen from Table~\ref{table4}, the new
method outperforms DIRECT and DIRECT{\it l} also on Criterion C3.

\begin{table}[t]
\caption{Number of hyperintervals for GKLS test functions
(Criterion C2).}
\begin{center} \scriptsize \label{table3}
\begin{tabular}{@{\extracolsep{\fill}}|c|c|c|r|r|r|r|r|r|}\hline
$N$ & $\Delta$ & Class &
\multicolumn{3}{c|}{50\%} & \multicolumn{3}{c|}{100\%}\\
\cline{4-9}& & &DIRECT &DIRECT{\it l} &New &DIRECT &DIRECT{\it l} &New\\
\hline
2 &$10^{-4}$ & simple &111  &152  &185 &1159 &2318 &1137\\
2 &$10^{-4}$ & hard  &1062 &1328 &607 &3201 &3414 &3993\\
\hline
3 &$10^{-6}$ & simple &386  &591  &1867 &12507 &13309 &12149 \\
3 &$10^{-6}$ & hard &1749 &1967 &2061 &$>$1000000 (4) &29233 &14357\\
\hline
4 &$10^{-6}$ & simple &4805  &7194  &21635 &$>$1000000 (4) &118744 &186295 \\
4 &$10^{-6}$ & hard &16114 &33147 &33173 &$>$1000000 (7) &287857 &223263\\
\hline
5 &$10^{-7}$ & simple &1660   &9246   &19823  &$>$1000000 (1)&178217     &255059 \\
5 &$10^{-7}$ & hard &55092  &126304 &169413 &$>$1000000 (16) &$>$1000000 (4) &1592969 \\
\hline
\end{tabular}
\end{center}
\end{table}

\begin{table}[!t]
\caption{Average number of trial points for GKLS test functions
(Criterion C3).}
\begin{center} \scriptsize \label{table4}
\begin{tabular}{@{\extracolsep{\fill}}|c|c|c|r|r|r|c|c|}\hline
$N$ & $\Delta$ & Class &DIRECT\hspace{1mm}
&DIRECT{\it l}\hspace{1mm} &New\hspace{2mm} &\multicolumn{2}{c|}{Improvement} \\
\cline{7-8} & & & & & &DIRECT/New &DIRECT{\it l}/New \\ \hline
2 &$10^{-4}$ & simple &198.89   &292.79  &97.22  &2.06 &3.01 \\
2 &$10^{-4}$ & hard &1063.78  &1267.07 &192.00 &5.54 &6.60 \\
\hline
3 &$10^{-6}$ & simple &1117.70       &1785.73  &491.28 &2.28     &3.63 \\
3 &$10^{-6}$ & hard & $>$42322.65  &4858.93  &618.32 &$>$68.45 &7.86 \\
\hline
4 &$10^{-6}$ & simple &$>$47282.89   &18983.55  &3675.84 &$>$12.87 &5.16 \\
4 &$10^{-6}$ & hard &$>$95708.25   &68754.02  &5524.77 &$>$17.32 &12.44 \\
\hline
5 &$10^{-7}$ & simple &$>$16057.46   &16758.44      &3759.05  &$>$4.27 &4.46     \\
5 &$10^{-7}$ & hard &$>$217215.58  &$>$269064.35  &22189.47 &$>$9.79 &$>$12.13 \\
\hline
\end{tabular}
\end{center}
\end{table}

Finally, results of comparison between the new algorithm and its
two competitors in terms of Criterion C4 are reported in
Table~\ref{table5}. This table shows how often the new algorithm
was able to minimize each of 100 functions of a class with a
smaller number of trials with respect to DIRECT or DIRECT{\it l}.
The notation `$p$\,:\,$q$' means that among~100 functions of a
particular test class there are $p$ functions for which DIRECT (or
DIRECT{\it l}) spent fewer function trials than the new algorithm
and $q$ functions for which the new algorithm generated fewer
trial points with respect to DIRECT (or DIRECT{\it l}) ($p$ and
$q$ are from~(\ref{winDIRECT}) and~(\ref{winNEW}), respectively).
As a rule, the more hard objective functions are presented in a
test class, the more pronounced becomes the advantage of the new
algorithm on Criterion C4, as well.

\begin{table}[t]
\caption{Comparison between the new algorithm and DIRECT and
DIRECT{\it l} in terms of Criterion C4.}
\begin{center}\scriptsize \label{table5}
\begin{tabular}{@{\extracolsep{\fill}}|c|c|c|c|c|}\hline
$N$ & $\Delta$ & Class &DIRECT\,:\,New &DIRECT{\it l}\,:\,New \\
\hline
2 &$10^{-4}$ & simple  & 28\,:\,72 & 21\,:\,79  \\
2 &$10^{-4}$ & hard  & 15\,:\,85 & 16\,:\,84  \\
\hline
3 &$10^{-6}$ & simple  & 36\,:\,64 & 30\,:\,70  \\
3 &$10^{-6}$ & hard  & 19\,:\,81 & 17\,:\,83  \\
\hline
4 &$10^{-6}$ & simple  & 39\,:\,61 & 25\,:\,75  \\
4 &$10^{-6}$ & hard  & 14\,:\,86 & 16\,:\,84  \\
\hline
5 &$10^{-7}$ & simple  & 55\,:\,45 & 17\,:\,83  \\
5 &$10^{-7}$ & hard  & 26\,:\,74 & 20\,:\,80  \\
\hline
\end{tabular}
\end{center}
\end{table}

As demonstrated by the results of the extensive numerical
experiments performed, the usage of the gradient information
together with the efficient partitioning strategy allows one to
obtain a serious acceleration in comparison with the DIRECT-based
methods on the studied classes of test problems.

\bibliographystyle{amsplain}

\begin{thebibliography}{10}

\bibitem{Audet:et:al.(2005)}
C.~Audet, P.~Hansen, and G.~Savard (eds.), \emph{Essays and surveys
in global
  optimization}, {GERAD} 25th Anniversary, Springer--Verlag, New York, 2005.

\bibitem{Baritompa(1993)}
W.~Baritompa, \emph{Customizing methods for global optimization --
{A}
  geometric viewpoint}, J. Global Optim. \textbf{3} (1993), no.~2, 193--212.

\bibitem{Breiman&Cutler(1993)}
L.~Breiman and A.~Cutler, \emph{A deterministic algorithm for global
  optimization}, Math. Program. \textbf{58} (1993), no.~1--3, 179--199.

\bibitem{DiSerafino:et:al.(2011)}
D.~{Di Serafino}, G.~Liuzzi, V.~Piccialli, F.~Riccio, and
G.~Toraldo, \emph{A
  modified {DIviding} {RECTangles} algorithm for a problem in astrophysics}, J.
  Optim. Theory Appl. \textbf{151} (2011), no.~1, 175--190.

\bibitem{Dumas:et:al.(2009)}
L.~Dumas, B.~Druez, and N.~Lecerf, \emph{A fully adaptive hybrid
optimization
  of aircraft engine blades}, J. Comput. Appl. Math. \textbf{232} (2009),
  no.~1, 54--60.

\bibitem{Evtushenko(1985)}
{Yu.}~G. Evtushenko, \emph{Numerical optimization techniques},
Translations
  Series in Mathematics and Engineering, Springer--Verlag, Berlin, 1985.

\bibitem{Evtushenko:et:al.(2009)}
{Yu.}~G. Evtushenko, V.~U. Malkova, and A.~A. Stanevichyus,
\emph{Parallel
  global optimization of functions of several variables}, Comput. Math. Math.
  Phys. \textbf{49} (2009), no.~2, 246--260.

\bibitem{Evtushenko&Posypkin(2011)}
{Yu.}~G. Evtushenko and M.~A. Posypkin, \emph{An application of the
nonuniform
  covering method to global optimization of mixed integer nonlinear problems},
  Comput. Math. Math. Phys. \textbf{51} (2011), no.~8, 1286--1298.

\bibitem{Finkel&Kelley(2006)}
D.~E. Finkel and C.~T. Kelley, \emph{Additive scaling and the
{DIRECT}
  algorithm}, J. Global Optim. \textbf{36} (2006), no.~4, 597--608.

\bibitem{Gablonsky&Kelley(2001)}
J.~M. Gablonsky and C.~T. Kelley, \emph{A locally-biased form of the
{DIRECT}
  algorithm}, J. Global Optim. \textbf{21} (2001), no.~1, 27--37.

\bibitem{Gaviano&Lera(2002)}
M.~Gaviano and D.~Lera, \emph{A complexity analysis of local search
algorithms
  in global optimization}, Optim. Methods Softw. \textbf{17} (2002), no.~1,
  113--127.

\bibitem{Gaviano&Lera(2008)}
\bysame, \emph{A global minimization algorithm for {Lipschitz}
functions},
  Optim. Lett. \textbf{2} (2008), no.~1, 1--13.

\bibitem{Gaviano:et:al.(2003)}
M.~Gaviano, D.~Lera, D.~E. Kvasov, and {Ya.}~D. Sergeyev,
\emph{Algorithm 829:
  Software for generation of classes of test functions with known local and
  global minima for global optimization}, ACM Trans. Math. Software \textbf{29}
  (2003), no.~4, 469--480.

\bibitem{Gergel(1997)}
V.~P. Gergel, \emph{A global optimization algorithm for multivariate
function
  with {Lipschitzian} first derivatives}, J. Global Optim. \textbf{10} (1997),
  no.~3, 257–--281.

\bibitem{Graf:et:al.(2007)}
P.~A. Graf, K.~Kim, W.~B. Jones, and {L.-W.} Wang, \emph{Surface
passivation
  optimization using {DIRECT}}, J. Comput. Phys. \textbf{224} (2007), no.~2,
  824--835.

\bibitem{He:et:al.(2002)}
J.~He, L.~T. Watson, N.~Ramakrishnan, C.~A. Shaffer, A.~Verstak,
J.~Jiang,
  K.~Bae, and W.~H. Tranter, \emph{Dynamic data structures for a direct search
  algorithm}, Comput. Optim. Appl. \textbf{23} (2002), no.~1, 5--25.

\bibitem{Horst&Pardalos(1995)}
R.~Horst and P.~M. Pardalos (eds.), \emph{Handbook of global
optimization},
  vol.~1, Kluwer Academic Publishers, Dordrecht, 1995.

\bibitem{Horst&Tuy(1996)}
R.~Horst and H.~Tuy, \emph{Global optimization -- deterministic
approaches},
  Springer--Verlag, Berlin, 1996.

\bibitem{Jones:et:al.(1993)}
D.~R. Jones, C.~D. Perttunen, and B.~E. Stuckman, \emph{Lipschitzian
  optimization without the {Lipschitz} constant}, J. Optim. Theory Appl.
  \textbf{79} (1993), no.~1, 157--181.

\bibitem{Jones:et:al.(1998)}
D.~R. Jones, M.~Schonlau, and W.~J. Welch, \emph{Efficient global
optimization
  of expensive black-box functions}, J. Global Optim. \textbf{13} (1998),
  no.~4, 455--492.

\bibitem{Kvasov:et:al.(2003)}
D.~E. Kvasov, C.~Pizzuti, and {Ya.}~D. Sergeyev, \emph{Local tuning
and
  partition strategies for diagonal {GO} methods}, Numer. Math. \textbf{94}
  (2003), no.~1, 93--106.

\bibitem{Kvasov&Sergeyev(2003)}
D.~E. Kvasov and {Ya.}~D. Sergeyev, \emph{Multidimensional global
optimization
  algorithm based on adaptive diagonal curves}, Comput. Math. Math. Phys.
  \textbf{43} (2003), no.~1, 42--59.

\bibitem{Kvasov&Sergeyev(2009)}
\bysame, \emph{A univariate global search working with a set of
{Lipschitz}
  constants for the first derivative}, Optim. Lett. \textbf{3} (2009), no.~2,
  303--318.

\bibitem{Lera&Sergeyev(2010a)}
D.~Lera and {Ya.}~D. Sergeyev, \emph{An information global
minimization
  algorithm using the local improvement technique}, J. Global Optim.
  \textbf{48} (2010), no.~1, 99--112.

\bibitem{Lera&Sergeyev(2010b)}
\bysame, \emph{Lipschitz and {H\"{o}lder} global optimization using
  space-filling curves}, Appl. Numer. Math. \textbf{60} (2010), no.~1--2,
  115--129.

\bibitem{Liuzzi:et:al.(2010)}
G.~Liuzzi, S.~Lucidi, and V.~Piccialli, \emph{A partition-based
global
  optimization algorithm}, J. Global Optim. \textbf{48} (2010), no.~1,
  113--128.

\bibitem{Luo:et:al.(2011)}
C.~Luo, {S.-L.} Zhang, C.~Wang, and Z.~Jiang, \emph{A
metamodel-assisted
  evolutionary algorithm for expensive optimization}, J. Comput. Appl. Math.
  \textbf{236} (2011), no.~5, 759--764.

\bibitem{Mockus(2000)}
J.~Mockus, \emph{A set of examples of global and discrete
optimization:
  {Applications} of bayesian heuristic approach}, Kluwer Academic Publishers,
  Dordrecht, 2000.

\bibitem{Moles:et:al.(2003)}
C.~G. Moles, P.~Mendes, and J.~R. Banga, \emph{Parameter estimation
in
  biochemical pathways: {A} comparison of global optimization methods}, Genome
  Res. \textbf{13} (2003), no.~11, 2467--2474.

\bibitem{Nesterov(2004)}
Yu. Nesterov, \emph{Introductory lectures on convex optimization: {A
Basic}
  course}, Kluwer Academic Publishers, Dordrecht, 2004.

\bibitem{Nocedal&Wright(1999)}
J.~Nocedal and S.~J. Wright, \emph{Numerical optimization},
Springer--Verlag,
  Dordrecht, 1999.

\bibitem{Panning:et:al.(2008)}
T.~D. Panning, L.~T. Watson, N.~A. Allen, K.~C. Chen, C.~A. Shaffer,
and J.~J.
  Tyson, \emph{Deterministic parallel global parameter estimation for a model
  of the budding yeast cell cycle}, J. Global Optim. \textbf{40} (2008), no.~4,
  719--738.

\bibitem{Pardalos:et:al.(2000)}
P.~M. Pardalos, H.~E. Romeijn, and H.~Tuy, \emph{Recent developments
and trends
  in global optimization}, J. Comput. Appl. Math. \textbf{124} (2000), no.~1-2,
  209--228.

\bibitem{Pinter(1996)}
J.~{Pint\'{e}r}, \emph{Global optimization in action (continuous and
lipschitz
  optimization: Algorithms, implementations and applications)}, Kluwer Academic
  Publishers, Dordrecht, 1996.

\bibitem{Sergeyev(1995)}
{Ya.}~D. Sergeyev, \emph{An information global optimization
algorithm with
  local tuning}, SIAM J. Optim. \textbf{5} (1995), no.~4, 858--870.

\bibitem{Sergeyev(1998)}
\bysame, \emph{Global one-dimensional optimization using smooth
auxiliary
  functions}, Math. Program. \textbf{81} (1998), no.~1, 127--146.

\bibitem{Sergeyev(2000)}
\bysame, \emph{An efficient strategy for adaptive partition of
  {$N$}-dimensional intervals in the framework of diagonal algorithms}, J.
  Optim. Theory Appl. \textbf{107} (2000), no.~1, 145--168.

\bibitem{Sergeyev(2005)}
\bysame, \emph{Efficient partition of {$N$}-dimensional intervals in
the
  framework of one-point-based algorithms}, J. Optim. Theory Appl. \textbf{124}
  (2005), no.~2, 503--510.

\bibitem{Sergeyev:et:al.(1999)}
{Ya.}~D. Sergeyev, P.~Daponte, D.~Grimaldi, and A.~Molinaro,
\emph{Two methods
  for solving optimization problems arising in electronic measurements and
  electrical engineering}, SIAM J. Optim. \textbf{10} (1999), no.~1, 1--21.

\bibitem{Sergeyev&Kvasov(2006)}
{Ya.}~D. Sergeyev and D.~E. Kvasov, \emph{Global search based on
efficient
  diagonal partitions and a set of {Lipschitz} constants}, SIAM J. Optim.
  \textbf{16} (2006), no.~3, 910--937.

\bibitem{Sergeyev&Kvasov(2008)}
\bysame, \emph{Diagonal global optimization methods}, {FizMatLit},
Moscow,
  2008, In Russian.

\bibitem{Sergeyev&Kvasov(2011)}
{Ya}.~D. Sergeyev and D.~E. Kvasov, \emph{Lipschitz global
optimization}, Wiley
  Encyclopedia of Operations Research and Management Science (J.~J. Cochran,
  ed.), vol.~4, Wiley, New York, 2011, pp.~2812--2828.

\bibitem{Strongin&Sergeyev(2000)}
R.~G. Strongin and {Ya.}~D. Sergeyev, \emph{Global optimization with
non-convex
  constraints: {Sequential} and parallel algorithms}, Kluwer Academic
  Publishers, Dordrecht, 2000.

\bibitem{Trigiante(2000)}
D.~Trigiante (ed.), \emph{Recent trends in numerical analysis}, Nova
Science
  Publishers, Inc., New York, 2000.

\bibitem{Wu:et:al.(2005)}
Y.~Wu, L.~Ozdamar, and A.~Kumar, \emph{{TRIOPT}: {A}
triangulation-based
  partitioning algorithm for global optimization}, J. Comput. Appl. Math.
  \textbf{177} (2005), no.~1, 35--53.

\bibitem{Zhigljavsky&Zilinskas(2008)}
A.~A. Zhigljavsky and A.~{\v{Z}ilinskas}, \emph{Stochastic global
  optimization}, Springer, New\,York, 2008.

\end{thebibliography}

\providecommand{\bysame}{\leavevmode\hbox
to3em{\hrulefill}\thinspace}
\providecommand{\MR}{\relax\ifhmode\unskip\space\fi MR }
\providecommand{\MRhref}[2]{%
  \href{http://www.ams.org/mathscinet-getitem?mr=#1}{#2}
} \providecommand{\href}[2]{#2}

\end{document}